\documentclass [10pt,reqno]{amsart}
\usepackage{ucs}
\usepackage{amsmath, amssymb, amscd}
\usepackage{pb-diagram}
\usepackage[latin1]{inputenc}
\usepackage{amsmath}
\usepackage{amssymb}

\newtheorem{theorem}{Theorem}[section]
\newtheorem{definition}[theorem]{Definition}
\newtheorem{lemma}[theorem]{Lemma}
\newtheorem{prop}[theorem]{Proposition}
\newtheorem{corollary}[theorem]{Corollary}

\newtheorem*{remark}{Remark}

\newcommand{\abs}[1]{\lvert#1\rvert}
\newcommand{\norm}[1]{\| #1 \|}
\newcommand{\klammer}[1]{\left( #1 \right)}
\newcommand{\eckklammer}[1]{\left[ #1 \right]}

\newcommand{\geklatau}[1]{\left\{ #1 \right\}_\tau}
\newcommand{\dbar}{\bar{\partial}}
\renewcommand{\epsilon}{\varepsilon}

\DeclareMathAlphabet{\mathpzc}{OT1}{pzc}{m}{it}
\usepackage{mathrsfs}

\newcommand{\Z}{\mathbb{Z}}
\newcommand{\C}{\mathbb{C}}

\newcommand{\R}{\mathbb{R}}

\renewcommand{\qed}{$\hfill \square$ \bigskip \\}
\renewcommand{\phi}{\varphi}

\newcommand{\ad}{\text{ad}}

\newcommand{\tensor}{\otimes}
\newcommand{\End}{\text{End}}
\newcommand{\Hom}{\text{Hom}}

\renewcommand{\det}{\text{det}}
\newcommand{\id}{\text{id}}

\newcommand{\su}{\mathfrak{su}}
\renewcommand{\sl}{\mathfrak{sl}}
\newcommand{\gl}{\mathfrak{gl}}
\renewcommand{\u}{\mathfrak{u}}

\renewcommand{\bar}{\overline}
\newcommand{\G}{\mathscr{G}}

\newcommand{\tr}{\text{tr}}
\newcommand{\data}{\mathfrak{s},E}

\newcommand{\conf}{\mathscr{C}}
\newcommand{\bonf}{\mathscr{B}}
\newcommand{\s}{\mathfrak{s}}

\newcommand{\hata}{{\hat{A}}}

\begin{document}

\thispagestyle{empty}

\title [What to expect from $U(n)$ Seiberg-Witten monopoles for $n > 1$]
{What to expect from $U(n)$ Seiberg-Witten monopoles for $n > 1$} \author
[Raphael
Zentner]{Raphael Zentner} 

\begin {abstract} 
We study generalisations to the structure groups U(n) of the familiar (abelian) Seiberg-Witten monopole equations on a four-manifold $X$ and their moduli spaces. For $n=1$ one
obtains the classical monopole equations. For $n > 1$ our results indicate that
there should not be any non-trivial gauge-theoretical invariants which are obtained by the scheme `evaluation of cohomology classes on the fundamental cycle of the moduli space'. For, if $b_2^+$ is positive the moduli space should be
`cobordant' to the empty space because we can deform the equations so as the moduli space of the deformed equations is generically empty. Furthermore, on K\"ahler surfaces
with $b_2^+ > 1$, the moduli spaces
become empty as soon as we perturb with a non-vanishing holomorphic 2-form.
\end {abstract}

\address {Fakult\"at f\"ur Mathematik \\ Universit\"at Bielefeld \\ 33501 Bielefeld\\
Germany}

\email{rzentner@math.uni-bielefeld.de}

\maketitle

\section*{Introduction}
In this paper we study generalisations of the familiar (abelian) Seiberg-Witten monopole equations to the structure groups $U(n)$ for $n > 1$. This is done by twisting a given $Spin^c$ structure $\mathfrak{s}$ on the four-manifold $X$ with a Hermitian bundle $E$ of rank $n$. The first variable in the theory will consist of sections $\Psi$ of the twisted spinor bundle $S^+_\mathfrak{s} \tensor E$. As for the second, 
 instead of taking the $Spin^c$ connections in $\mathfrak{s}$ as variables in the theory, we keep one fixed as a parameter and take the $U(n)$ connections $\hata$ in the bundle $E$ as a variable. There are then straightforward generalisations of the classical Seiberg-Witten equations to this situation:
\begin{equation*}
\begin{split}
 D_{\hat{A}} \Psi & =  0 \\
 \gamma(F_{\hat{A}}^+) - \mu_{0,\tau}(\Psi) & =  \gamma(\eta) \ \id   \ .
\end{split}
\end{equation*}
Here $D_{\hata}$ is the associated Dirac-operator to $\hata$, the map $\gamma$ is derived from Clifford-multiplication, $F_\hata^+$ is the self-dual part of the curvature of the connection $\hata$, $\eta$ is a self-dual 2-form which serves as a perturbation of the equations, and $\mu_{0,\tau}$ is a quadratic map in the spinor depending on a parameter $\tau \in [0,1]$, explicitely described below. 

The involved analysis is much more difficult than in classical Seiberg-Witten theory. First, the associated moduli spaces are in general not compact anymore but have a natural compactification similar to the Uhlenbeck-compactification of the moduli spaces of instantons. Second, generic regularity is a much harder problem than in classical Seiberg-Witten theory and cannot be achieved by the perturbation with the self-dual 2-form $\eta$ above (or the metric in addition). Third, we cannot avoid `reducibles' in general, i.e. solutions $(\Psi,\hata)$ to the above equations which have positive-dimensional stabiliser under the action of the gauge group. 

The aim of the present paper is easily stated: Without even solving all of the mentioned technical problems we shall show that it is not really worth to do so, because we get quite strong evidence that there should not be interesting gauge-theoretical invariants envolved, at least none which are derived with the classical scheme `evaluation of cohomology classes on the fundamental cycle of the moduli space'. This evidence is given by two main results. The first, Proposition \ref{main result 1} states that if we put $\tau = 0$ in the above equation then for a generic perturbation $\eta$ the associated moduli space is empty if $b_2^+(X) \geq 1$. It should be pointed out that putting $\tau = 0$ is only sensible for $n > 1$ because otherwise we lose control over the compactness or compactification. But moduli spaces for different $\tau$ should be `cobordant' if an invariant is defined at all.
The second, Corollary \ref{final conclusion}, shows that on a K\"ahler surface the moduli space is empty as soon as one perturbs with a non-vanishing holomorphic 2-form, which is always possible if $b_2^+(X) > 1$.

Our reason for studying $U(n)$ monopoles consists in the fact that these appeared naturally when studying certain $PU(N)$ monopoles, for integers $1 \leq n < N$. The $PU(2)$ monopoles have been used extensively with the aim of proving Witten's conjecture \cite{W} on the relation between the Seiberg-Witten and the Donaldson invariants, first by Pidstrigach-Tyurin \cite{PT}, Okonek-Teleman \cite{OT}, \cite{T2}, and then by Feehan and Leness \cite{FL1}, \cite{FL2}, \cite{FL3}, \cite{FL4}. Feehan and Leness now seem to have proved the full conjecture \cite{FL5}. 

Kronheimer has introduced instanton-type invariants associated to Hermitian bundles of rank $N$ \cite{K} which are a generalisation of the polynomial invariants of Donaldson appearing when $N=2$. Before these invariants were even properly defined the physicists Mari\~no and Moore conjectured that such invariants should not contain new differential topological information and suggested a generalisation of Witten's conjecture to a relationship between these invariants and the Seiberg-Witten invariants. Kronheimer verified this conjecture for a large class of four-manifolds. We, instead, have investigated a generalisation of the above mentioned approach by means of $PU(N)$ monopoles  \cite{Z2}. The main results in this paper, Proposition \ref{main result 1} and Corollary \ref{final conclusion} will be used in \cite{Z2} in order to give first steps towards a proof of the mentioned conjecture. Aside this motivation, studying $U(n)$ monopoles is also interesting in itself.

In the first section we shall introduce our setting, define the above mentioned map $\mu_{0,\tau}$ and derive some important properness property which will imply the existence of an a-priori $C^0$ bound on the spinor component of a monopole. Given this bound one can show that there is a natural Uhlenbeck-type compactification of the moduli space. In the second section we discuss the implications of deforming the equations by $\tau \in [0,1]$, yielding the above mentioned first main result. In the third section we discuss $U(n)$ moduli spaces on the K\"ahler surfaces yielding the mentioned vanishing result.
 \\

\section*{Acknowledgements}
I am grateful to Andrei Teleman for many mathematical discussions. This work owes much to him. I would also like to thank Peter Kronheimer and Kim Fr\o yshov for discussions on some of the rather technical aspects of the theory.

\section{The $U(n)$-monopole equations}
In this section we shall define the $U(n)$ monopole equations, study some of their basic properties and define the moduli space. The standard material in Seiberg-Witten theory ($Spin^c$ structures, $Spin^c$ connections etc.) can be found in one of the textbooks on the topics like \cite{N}, \cite{M} or diverse lecture notes like \cite{T3}.

\subsection{The configuration space}
Let $X$ be a closed oriented Riemannian four-manifold with a
$Spin^c$ structure $\mathfrak{s}$ on it. The $Spin^c$ structure consists of two
Hermitian rank 2 vector bundles $S^\pm_{\mathfrak{s}}$ with identified
determinant line bundles and a Clifford multiplication 
\begin{equation*} 
\gamma : \Lambda^1(T^*X) \to \Hom_\C(S^+_\mathfrak{s},S^-_\mathfrak{s}) \ .
\end{equation*}

Furthermore suppose we are given a
Hermitian vector bundle $E$ with determinant line bundle $w = \det(E)$ on $X$.
We can then form spinor bundles 
\begin{equation*}
W_{\data}^\pm := S^\pm_{\mathfrak{s}} \tensor E .
\end{equation*}
Clifford multiplication extends by tensoring with the identity on $E$. 
This way we obtain a $Spin^c$ - structure `twisted' by the hermitian bundle $E$.
\\

We shall denote by $\mathscr{A}(E)$ the space of smooth unitary connections on
$E$ which is an affine space modelled on
$\Omega^1(X;\mathfrak{u}(E))$. Here $\u(E)$ denotes the bundle of skew-adjoint
endomorphisms of $E$. Furthermore $\Gamma(X;W^+_{\data})$ denotes the space of
smooth sections of the spinor bundle $W^+_{\data}$.
We define our configuration space to be the space
\begin{equation*}
\mathscr{C}_{\data} := \Gamma(X;W^+_{\data}) \times \mathscr{A}(E) \ .
\end{equation*}
We denote by  $\G$ the group of
unitary automorphisms of $E$; it is the `gauge group' of our problem. It acts
 in a canoncial way on sections of the spinor bundles, and as $(u,\nabla_A)
\mapsto u \nabla_A u^{-1}$ on the connections, where $u$ is a gauge
transformation and $\nabla_A$ a unitary connection. The set
$\mathscr{B}_{\data}$ is defined to be the configuration space up to gauge, i.e. the quotient space $\conf_{\data} / \G$. 
\\

The reason we consider only smooth objects is purely a matter of simplicity
here. Obviously, as soon as we wish to consider more analytical properties like
transverality, we study suitable Sobolev-completions of these
spaces.
\\

A configuration $(\Psi,\hat{A})$ shall be called irreducible if its stabiliser
$\Gamma(\Psi,\hat{A}) \subseteq \G $ inside the gauge group is trivial. We
shall denote by $\conf^*_{\data}$ the open subspace of configurations with
trivial stabiliser and by $\bonf^*_{\data}$ its quotient space. If we consider
 Sobolev completions of our spaces it is standard
to show that there are local slices for the $\G$ action on $\conf^*_{\data}$, so
that $\bonf^*_{\data}$ becomes a Banach manifold.
\\

\subsection{Algebraic preliminaries}

We shall now define a quadratic map in the spinor that will
appear in the $U(n)$ monopole equations. We prove a properness property for
this map that will be essential in proving a uniform bound for solutions to the
monopole equations. \\

The twisted spinor bundles
$W^\pm_{\data}$ are associated bundles of the fibre product of a $Spin^c$
principal bundle and a $U(n)$-principal bundle on X, with the standard fibre
$\C^2 \tensor \C^n$.
Let us consider the isomorphism
\begin{equation*}
\begin{split}
  (p,q): \ \,  \gl(\C^n) & \to \sl(\C^n) \oplus \C \, \id \\
   	 a  & \mapsto  \left(a - \frac{1}{n} \ \tr(a) \cdot \id, \frac{1}{n} \
\tr(a) \cdot \id\right) \ .
\end{split}
\end{equation*}
Both components $p$ and $q$ are orthogonal projections onto their images. Note
that $\gl(\C^2) \tensor \gl(\C^n)$ and $\gl(\C^2 \tensor \C^n)$ are canonically
isomorphic. We define the orthogonal projections
\begin{equation*}
\begin{split}
  P: \ \gl(\C^2 \tensor \C^n) & \to \sl(\C^2) \tensor \sl(\C^n) \ , \\
  Q: \ \gl(\C^2 \tensor \C^n) & \to \sl(\C^2) \tensor \C \, \id
\end{split}
\end{equation*}
to be the tensor product $( \ )_0 \tensor p$ respectively $( \ )_0 \tensor q$,
with $( \ )_0$ denoting the trace-free part of the endomorphism of the first
factor $\C^2$. 
 
For elements $\Psi, \Phi \in \C^2 \tensor \C^n$ we define the map
\[
\mu_{0,\tau}(\Psi,\Phi):=P (\Psi \Phi^*) \ + \tau \, Q ( \Psi \Phi^*) ,
\]
 where $(\Psi \Phi^*) \in \gl(\C^2 \tensor \C^n)$ is defined to be the
endomorphism $\Xi \mapsto
(\Phi,\Xi) \cdot \Psi$. Furthermore $\tau  \in [0,1]$.

With this notation $\mu_{0,1} (\Psi,\Phi)$ is simply the orthogonal projection
of the endomorphism $\Psi \Phi^* \in \gl(\C^2 \tensor \C^n)$ onto $\sl(\C^2)
\tensor \gl(\C^n)$. We shall also write $\mu_{0,\tau}(\Psi):=
\mu_{0,\tau}(\Psi,\Psi)$ for the associated quadratic map. In the case $n=1$ the
map $\mu_{0,1}(\Psi)$ is the quadratic map in the spinor usually occuring in the
Seiberg-Witten equations \cite{K}. \\
%

\begin{prop}\label{properness}
Suppose $n > 1$. Then the quadratic map $\mu_{0,\tau}$ is uniformly proper. In
other words, there is a positive
constant $c > 0$ such that 
  \begin{equation}\label{propernessconstant}
    \abs{\mu_{0,\tau}(\Psi)} \geq c \abs{\Psi}^2 \ , 
  \end{equation}
independently of $\tau \in [0,1]$. 
As a consequence we have the formula
  \begin{equation}
    \left(\mu_{0,\tau}(\Psi) \Psi, \Psi \right) \geq c^2 \abs{\Psi}^4 \ 
  \end{equation}
whenever $\tau \geq 0$.

\end{prop}
We defer the proof of this proposition to the appendix. \qed 

Because of the equivariance property of the map $\mu_{0,\tau}$ we get in
a straightforward way corresponding maps between bundles, giving rise to
\begin{equation*}
\mu_{0,\tau}: W^\pm_{\data}  \times  W^\pm_{\data} \ 
\to \sl(S^\pm_{\s}) \tensor_{\C} \gl(E) \ ,
\end{equation*}
respectively, for the quadratic map,
\begin{equation*}
\mu_{0,\tau}: W^\pm_{\data} \ 
\to \su(S^\pm_{\s}) \tensor_{\R} \u(E) \ .
\end{equation*}
These maps on the bundle level satisfy the corresponding statement in the above
proposition with the same constant $c$. 

\subsection{The $U(n)$-monopole equations}

The Clifford map $\gamma$ is, up to a universal constant, an isometry of the
cotangent bundle
onto a real form inside $\Hom_\C(S^+_\mathfrak{s},S^-_\mathfrak{s})$ which can
be specified by the Pauli matrices. We extend $\gamma$ to $\End(S^+_\s
\oplus S^-_\s)$ by $ - \gamma^*$ on the negative Spinor bundle. It then
naturally extends to exteriour powers of $T^*X$, and in particular its
restriction to self-dual two-forms is zero on the negative Spinor bundle, and
induces an isomorphism 
\begin{equation*}
\gamma: \Lambda^2_+(T^*X) \stackrel{\cong}{\to} \su(S^+_\mathfrak{s}) \ .
\end{equation*}

Let's fix a background $Spin^c$ connection $B$ on $\mathfrak{s}$. By composing
the connection $\nabla_{B} \tensor \nabla_{\hat{A}}$ with the Clifford
multiplication we get a Dirac operator
\begin{equation*}
D_{\hat{A}} := \gamma \circ \left(\nabla_{B} \tensor \nabla_{\hat{A}}\right) :
\Gamma(X;W^\pm_{\data}) \to \Gamma(X;W^\mp_{\data}) \ .
\end{equation*}
Its extension to sections of $W^+_{\data} \oplus W^-_{\data}$  is a self-adjoint first order elliptic operator. We have
oppressed the $Spin^c$ connection $B$ from the notation because it will not be
a variable in our theory.

For a configuration $(\Psi,\hat{A}) \in \conf_{\data}$ the $U(n)$-monopole
equations with parameter $\tau \in [0,1]$ and perturbation $\eta \in \Omega^2_+(X;i\R)$
read
\begin{equation}\label{monopole equations}
\begin{split}
 D_{\hat{A}} \Psi & =  0 \\
 \gamma(F_{\hat{A}}^+) - \mu_{0,\tau}(\Psi) & =  \gamma(\eta) \ \id   \ .
\end{split}
\end{equation}
Here $F_{\hat{A}}$ designs the curvature of the connection $\hat{A}$ and
$F^+_{\hat{A}}$ its selfdual part. 

\subsection{The moduli space}
The left hand side of the above equations can be seen as a map $\mathscr{F}_\tau$
from the configuration space $\mathscr{C}_{\data}$ to
the space
$\Gamma(X;W^-_{\data})\times \Gamma(X;\su(S^+_{\mathfrak{s}}) \tensor
\u(E)) . $
This map satisfies
the equivariance property
\[
  \mathscr{F}_\tau(u.(\Psi,\hat{A})) = (u \times \ad_{u}) (\mathscr{F}_\tau(\Psi,\hat{A}))
\ 
\]
for $u \in \G$ and $(\Psi,\hat{A}) \in \conf_{\data}$. In particular, the set
of solutions to the above equations is gauge-invariant. The moduli space is
then defined to be the space of solutions to the monopole equations modulo
gauge:
\begin{equation*}
 M_{\data}(\tau,\eta) := \{[\Psi,\hat{A}] \in \bonf_{\data}
| \mathscr{F}_\tau(\Psi,\hat{A}) = (0,\gamma(\eta)) \} \ .
\end{equation*}
There is an elliptic deformation complex associated to a solution $x=(\Psi,\hat{A})$ of
the monopole equations. Let us denote by $\mathscr{C}_x^0 = \Gamma(X;\u(E))$
the Lie algebra of the gauge group, by $\mathscr{C}_x^1 = \Gamma(X; W^+_{\data}
\oplus\Lambda^1(X) \tensor \u(E))$ the tangent space to the configuration space
at $x$, and by $\mathscr{C}_x^2 = \Gamma(X;W^-_{\data} \oplus \Lambda^2_+(X)
\tensor \u(E))$ the target vector space of the monopole map $\mathscr{F}$.
Deriving the map $u \mapsto \mathscr{F}_\tau(u(\Psi,\hat{A}))$ yields then an 
elliptic deformation complex
\begin{equation*}
\begin{diagram}
\node{0} \arrow{e} \node{\mathscr{C}^0_x} \arrow{e,t}{\lambda_x}
\node{\conf^1_x} \arrow{e,t}{d_x \mathscr{F}}
\node{\conf^2_x} \arrow{e} \node{0 \ .}
\end{diagram}
\end{equation*}
Here $\lambda_x$ is the derivative of the map $u \mapsto u(\Psi,\hat{A})$, and
$d_x\mathscr{F}$ is the derivative of the monopole map $\mathscr{F}$ at the
solution $x=(\Psi,\hat{A})$. Let us denote by $H^i_x$ the associated cohomology
groups.
If the configuration $(\Psi,\hat{A})$ is
irreducible, that is, has trivial stabiliser, than the cohomology group $H^0_x$
vanishes. A solution is called {\em regular} if the second cohomology group
$H^2_x$ vanishes. If a solution $(\Psi,\hat{A})$ is regular and irreducible
then the local Kuranishi models for the moduli space show that the moduli space
is a smooth manifold in a neighbourhood of $[\Psi,\hat{A}]$, of dimension equal
to minus the index of the above elliptic complex. This `expected dimension' of the moduli space is computed with the Atiyah-Singer index theorem and is given by the following formula:
\begin{equation*}
\begin{split}
d(\data) := & -2 \,  \langle p_1(\su(E)),[X]\rangle - \,  n^2 (b_2^+(X) - b_1(X) + 1) \\
			& -\frac{n}{4} \, \text{sign}(X) + \langle \, c_1(E)^2 - 2\, c_2(E) + \, c_1(L)c_1(E) + \, \frac{n}{4} c_1(L)^2 , [X]\rangle \ ,
\end{split}
\end{equation*}
where $p_1(\su(E))$ denotes the first Pontryagin class of the bundle $\su(E)$, and $\text{sign}(X)$ the signature of the intersection form on $X$. In this formula the expression in the first line of the right hand side is the expected dimension of the moduli space of $U(N)$ ASD-connections in $E$, and the second line is the index of the Dirac operator $D_{\hata}$. 
\\

In general the moduli space $M_{\data}(\tau, \eta)$ will not consist only of irreducible
and regular solutions. One usually uses perturbations to the equations (or the 
map $\mathscr{F}$) in order to get a moduli space consisting of regular elements
only. For getting a well-defined moduli problem the perturbations need to
be equivariant as well. This makes `generic regularity' a harder problem in the
case $n > 1$ than
in the abelian situation $n=1$. 
The holonomy perturbations as appearing
in \cite{K} can be slightly modified to fit to our situation. It can then be  shown that for a generic perturbation the moduli space is a smooth manifold in
neighbourhoods of points $[\Psi,\hat{A}]$ for which we have $\Psi \neq 0$ and
$\hat{A}$ is an irreducible connection. In the instaton situation \cite{K} reducible connections can be generically avoided under suitable topological assumptions on $E$. This, however, does not seem to hold in the monopole situation in general. 

\subsection{Uniform bound on the spinor} 

For solutions to the monopole equations with $n \geq 2$ we will now deduce a
uniform bound on the spinor, which can be taken independently of the parameter
$\tau \geq 0$. For this, notice that the Weitzenb\"ock formula for the
Dirac-operator $D_{\hat{A}}$ reads
\begin{equation*}
\begin{split}
  D_{\hat{A}} D_{\hat{A}} & = \nabla_{B,\hat{A}}^* \nabla_{B,\hat{A}} +
\frac{1}{2} \gamma(F_{B,\hat{A}}) \\
			  & = \nabla_{B,\hat{A}}^* \nabla_{B,\hat{A}} +
\frac{s}{4} + \frac{1}{2} \gamma(\tr F_B) + \frac{1}{2} \gamma(F_{\hat{A}}) \ ,
\end{split}
\end{equation*}
where $s$ denotes the scalar curvature of the Riemannian four-manifold $X$,
$\nabla_{B,\hat{A}}$ the tensor product connection of the fixed $Spin^c$
connection $B$ and the $U(N)$ connection $\hat{A}$, and $F_{B,\hat{A}}$ its
curvature. Now suppose that we have a monopole $[\Psi,\hat{A}] \in
M_{\data}(\tau,\eta)$. Using the Weitzenb\"ock formula, the monopole
equation and the inequalities in the above proposition \ref{properness} now
yields the following inequality 
\begin{equation*} 
	\begin{split}
	\frac{1}{2} \Delta \abs{\Psi}^2 & = 
		\left( \nabla_{B,\hat{A}}^* \nabla_{B,\hat{A}}	\Psi, \Psi
		\right) \ - \ \abs{\nabla_{B,\hat{A}} \Psi}^2 \\
		& \leq -\frac{s}{4} \abs{\Psi}^2 \ 
			- \frac{1}{2} \left(\gamma(\tr F_B^+) \Psi, \Psi
				\right) \
			- \frac{1}{2} \left( \mu_{0,\tau}(\Psi) \Psi,\Psi
				\right) \
			- \frac{1}{2} \left( \gamma(\eta) \Psi,\Psi\right)\\
		& \leq \left( - \frac{s}{4} \ + \frac{1}{2} \abs{\tr F_B^+} 
			\ + \abs{\eta} \right) \abs{\Psi}^2 
			\ - \frac{c^2}{2}	\abs{\Psi}^4 \ .
	\end{split}
\end{equation*}
Let $K$ be the maximum over $X$ of the coefficient of $\abs{\Psi}^2$ in the last line. This quantity can be a priori negative.
At a point $x$ on the four-manifold where $\abs{\Psi}$ admits its maximum the
Laplacian $\Delta \abs{\Psi}^2$ must be positive. If $\abs{\Psi}^2(x) \neq 0 $
we may devide the above inequality by $\abs{\Psi}^2(x)$, yielding the desired
uniform bound:
\begin{prop}\label{uniform bound}
There are constants $c, K \in \R$, with $c > 0$, such that for any
monopole $[\Psi,\hat{A}]\in M_{\data}(\tau,\eta)$, with $\tau \geq 0$, we have
a $C^0$ bound:
\begin{equation*}
	\max \abs{\Psi}^2 \, \leq \, \max \left\{ 0, K/c^2 \right\} \ . 
\end{equation*}
Here the constant $K$ depends on the Riemannian metric, the fixed background
$Spin^c$ connection, and the perturbation form $\eta$ only, whereas the
constant $c$ is universal. In particular the bound is uniform in $\tau \in [0,1]$. If the constant $K$ is negative then there are only solutions with vanishing spinor component to the $U(N)$ monopole equations. 
\end{prop}

\subsection{Compactness}
Contrary to the Abelian case $n=1$ the $U(n)$ moduli spaces $M_{\data}(\tau,\eta)$ are in general not
compact. However, there is a natural  compactification of these moduli spaces similar to the
Uhlenbeck-compactification of instanton moduli spaces \cite{DK}. 
This subject has been treated with in detail in \cite{T2}, \cite{FL4} in the case of $PU(2)-$ monopoles. The
main reason why the Uhlenbeck-compactification carries over to the monopole situation is the uniform bound
on the spinor which we have dealt with in Proposition \ref{uniform bound} above. We will only describe this
compactification here and refer to the above mentioned references for the highly technical proofs. An outline of the proof in the $PU(n)$ situation can also be found in \cite{Z1}.

Let $\mathfrak{s}$ be a $Spin^c$-structure on $X$ and let $E \to X$ be a
unitary bundle on $X$. We denote by $E_{-k}$ a bundle which has first Chern
class $c_1(E_{-k}) = c_1(E)$ and whose second Chern class satisfies 
\[
\langle c_2(E_{-k}), [X] \rangle = \langle c_2(E), [X] \rangle - k \ .
\]
Such a bundle is unique up to isomorphism on a four-manifold.

\begin{definition}
An ideal monopole associated to the data $(\data)$ is given by a pair
$([\Psi,\hat{A}], \bf{x})$, where $[\Psi,\hat{A}] \in
M_{\mathfrak{s},E_{-k}}(\tau,\eta)$ is a $(\mathfrak{s},E_{-k})-$ monopole, and
$\bf{x}$ is an element in the k-th symmetric power $Sym^k(X)$ of X (that is, an
unordered set of k points in X, ${\bf x} = \{ x_1, \dots, x_k\}$). The
curvature density of $([\Psi,\hat{A}],{\bf x})$ is defined to be the measure
\[
\abs{F_\hata}^2 + 8 \pi^2 \sum_{x_i \in {\bf x}} \delta_{x_i} \ .
\]
The set of ideal monopoles associated to the data $(\data)$ and parameters $(\tau,\eta)$ is
\begin{equation}\label{idealmonopoles}
I M_{\data}(\tau,\eta) := \coprod_{k \geq 0} M_{\mathfrak{s},E_{-k}}(\tau,\eta)
\times Sym^k(X)
\end{equation}
\end{definition}

The set of ideal monopoles is then endowed with a convenient
topology. This is possible by specifying the underlying notion of convergence.
In fact, it can even be shown (cf. for instance \cite{T2}, p. 433, \cite{DK})
that this
topology can be induced by a metric on the set of ideal monopoles (but this
metric is not an extension of the natural metric induced by the $L^2$-metric on
the slices of
the gauge-group on the main-stratum). Either way, in this topology each
stratum has its natural topology, and we have the
following notion of
convergence of a sequence in the main-stratum:
\begin{definition}
The sequence of monopoles $[\Psi_n,\hat{A}_n] \in
\mathscr{M}_{\data}$ converges to the ideal monopole
$([\Psi,\hat{A}],{\bf x}) \in M_{\data_{-k}}(\tau,\eta) \times Sym^k(X)$ if we have:
\begin{enumerate}
\item The sequence of measures $\abs{F_{\hat{A}_n}}^2 \, vol_g $ converges as
measure to 
\[
  \abs{F_\hata}^2 \, vol_g + 8 \pi^2 \sum_{x_i \in {\bf x}} \delta_{x_i}
\]
\item On the complement $\Omega:= X - \{x \, | \, x \in {\bf x}\}$ there are
bundle
isomorphisms $u_{n} : E_{-k} |_{\Omega} \to E|_{\Omega}$ such that the
sequence 
\[
  u_{n}^* ( (\Psi_{n},\hat{A}_{n})|_{\Omega})
\]
converges in the $C^\infty$-topology to $(\Psi,\hat{A})|_{\Omega}$.
\end{enumerate}
\end{definition}

In a similar way the convergence of sequences in the lower strata are defined.

\begin{theorem}\label{compactness}
Let $[\Psi_n,\hat{A}_n] \in M_{\data}(\tau,\eta)$ be a sequence of
$U(N)$ monopoles. Then there is an integer $k\geq 0$, a multiset ${\bf x} \in
Sym^k(X)$ and a $U(N)$ monopole $[\Psi',\hat{A}'] \in M_{\data_{-k}}(\tau,\eta)$
such that the following is true: There is a subsequence $(n_k)$ such that the
sequence $[\Psi_{n_k},\hat{A}_{n_k}]$ converges to the ideal monopole
$([\Psi',\hat{A}'],{\bf x}) \in M_{\data_{-k}}(\tau,\eta) \times Sym^k(X)$ in
the above sense.
\end{theorem}

\begin{corollary}(Compactness-Theorem)
The closure of $M_{\data}$ inside the space of ideal monopoles
$IM_{\data} (\tau,\eta)$ is
compact. In fact, with the topology  specified on $IM_{\data}(\tau,\eta)$ this space is itself compact.
\end{corollary}

Suppose that we have a monopole $[\Psi,\hata] \in M_{\mathfrak{s},E}(\tau,\eta)$. We will show that the
$L^2-$norm of $F_\hata^+$ is uniformly bounded independently of the topological data $(\mathfrak{s},E)$. For
$\tau \in [0,1]$ we obviously have $\abs{\mu_{0,\tau}(\Psi)} \leq C \abs{\Psi}^2$ for a universal positive
constant $C$ (depending on $n$ only). From the curvature equation of the $U(n)-$monopole equations
(\ref{monopole equations}) we therefore get, together with the uniform bound (\ref{uniform bound}), a
pointwise inequality which integrated yields
\begin{equation}\label{bound Fplus}
\norm{F_\hata^+}^2_{L^2(X)} \leq \klammer{\frac{C K}{2 c^2} + \norm{\eta}_\infty}^2 vol_g(X) \ .
\end{equation}
From this and the Chern-Weil formulae it follows that the $L^2-$norm of the total curvature $F_\hata$ is
also bounded. This fact is an essential input for the proof of the compactness theorem.

\begin{remark}\label{cobordant}
We might also consider moduli spaces $\widetilde{M}_{\data} (\eta)$ of `pa\-ra\-met\-rised' $U(n)-$
monopoles, where as parameter we take $\tau \in [0,1]$. This parametrised moduli space then fibers over the
interval $[0,1]$. The fact that the uniform bound in (\ref{uniform bound}) can be taken independently of
$\tau$ indicates that we can also compactify the parametrised moduli space, and that fibrewise the
compactification coincides with the one considered above. Thus, heuristically, the compactified moduli
spaces $\overline{M}_{\data}(\tau,\eta)$ and $\overline{M}_{\data}(\tau',\eta)$ for $\tau,\tau' \in [0,1]$
are `cobordant'. 
\end{remark}


\section{A deformation of the equations for $n > 1$}
A natural idea is to study the dependence of the moduli space on the fixed parameter $\tau \in [0,1]$.
Surprisingly, we have the following result for the case $\tau = 0$:
\begin{prop}\label{main result 1}
Suppose the 4-manifold $X$ has $b_2^+(X)$ non-zero. Then for a generic
imaginary-valued self-dual 2-form $\eta$ the deformed moduli space
$M_{\data}(0,\eta)$ is empty.
\end{prop}
{\em Proof:} 
	Suppose $[\Psi,\hat{A}]$ belongs to the moduli space
$M_{\data}(0,\eta)$. In particular, the configuration $(\Psi,\hat{A})$ solves
the $U(n)$ - monopole equations (\ref{monopole equations}) with parameter
$\tau = 0$. Let us take the trace of the curvature equation in
(\ref{monopole equations}). We get, after applying the isomorphism $\gamma^{-1}$
the following formula:
	\begin{equation}\label{u1 asd}
		\tr (F_{\hat{A}}^+) = n \cdot \eta \ .
	\end{equation}
	But $\tr (F_{\hat{A}})$ is precisely the curvature $F_{\det(\hat{A})}$
of the connection $\det(\hat{A})$ that $\hat{A}$ induces in the determinant line
bundle $\det(E)$. Therefore the equation (\ref{u1 asd}) can be seen as a
perturbed $ASD$ - equation for connections in a line bundle. Now, if $A_0$ is a
fixed connection in the line bundle $\det(E)$, then any other connection $A$ is
given by $A_0 + a$, where $a$ is an imaginary valued one-form. Its curvature is
given by $F_A = F_{A_0} + da$, hence the equation $F_A^+ = n \cdot \eta$ is
equivalent to 
	\[
		d^+a = -F_{A_0}^+ + n \cdot \eta \ .
	\]
	But $d^+ : \Omega^1(X;i\R) \to \Omega^2_+(X;i\R)$ has cokernel
isomorphic to the space of self-dual harmonic imaginary-valued 2-forms, which is
of dimension $b_2^+(X)$. Thus, for generic $\eta \in \Omega^2_+(X;i\R)$ this
equation has no solution.
\qed
\begin{remark}
The Seiberg-Witten and Donaldson invariants are obtained by evaluating canonical cohomology classes on the
`fundamental cycle' given by the moduli space. Moduli spaces associated to different perturbations prove to
be cobordant under the condition $b_2^+(X) > 1$ and the canonical cohomology classes extend to the
cobordism.  The above proposition and remark \ref{cobordant} suggests that no non-trivial invariants of that
type should be expected from the $U(n)$ moduli spaces $M_{\data}(\tau,\eta)$.
\end{remark}

\section{$U(n)$ moduli spaces on K\"ahler surfaces}
In classical Seiberg-Witten theory K\"ahler surfaces are of a significant importance. Indeed, they provided the first examples of 4-manifolds with non-trivial Seiberg-Witten invariants \cite{W}. This was generalised to symplectic manifolds \cite{Ta}. All other non-vanishing results known to the author are derived from these manifolds by various kinds of glueing results for the Seiberg-Witten invariants \cite{Ta}, \cite{Fr}.

As the $U(n)$ monopole equations are a generalisation of the classical Seiberg-Witten equations it is therefore most natural to study the $U(n)$ monopole moduli spaces for K\"ahler surfaces. Whereas the analysis of the $U(n)$ monopole equations on K\"ahler surfaces is very analogous to the classical situation the final conclusion is in sharp contrast to the classical situation. Indeed, we will show in Corollary \ref{final conclusion} that if we perturb the monopole moduli space on a K\"ahler surface with a non-vanishing holomorphic 2-form then the associated moduli space is empty. 

Non-abelian monopoles on K\"ahler surfaces have also been studied by Teleman \cite{T}, Okonek and Teleman \cite{OT3} and by Bradlow and Garcia-Prada \cite{BG}, but with a rather complex geometric motivation. Corollary \ref{final conclusion} seems to appear here for the first time.

\subsection{The $U(n)$ - monopole equations on K\"ahler surfaces}
We will quickly recall now the canonical $Spin^c-$ structure on an almost complex surface. The additional
condition of $X$ being K\"ahler implies that there is a canonical $Spin^c$ connection induced by the
Levi-Civita connection. This will be our fixed back-ground $Spin^c$ connection and it is then simple to
determine the Dirac-operator associated to this fixed connection and a $U(n)-$ connection in a Hermitian
bundle $E$. We will then write down the $U(n)$ monopole equations in this particular setting.
\\

Suppose we have an almost complex structure $J: TX \to TX$ on the closed,
oriented Riemannian 4-manifold $X$ which is isometric. The associated K\"ahler
form
$\omega$ is defined by the formula
\[
   \omega_g(v,w):= g(Jv,w) \ .
\]
This is an anti-symmetric form of type $(1,1)$ when extended to the
complexification $TX^\C := TX
\tensor_\R \C$. It is a fundamental fact that the complexification of the bundle
of self-dual two forms is given by 
\[
  \Lambda^2_+ \tensor \C = \C \omega_g \oplus \Lambda^{2,0} \oplus \Lambda^{0,2}
\ .
\]
Let $e(u)$ denotes exterior multiplication with the form $u \in
\Lambda(T^*X^\C)$
and $e^*(u)$ its adjoint with respect to the inner product induced by the
Riemannian metric.

There is a canonical $Spin^c$-structure 
associated to an
almost-complex structure $J$ on $X$ \cite{Hi}. We shall denote it by
$\mathfrak{c}$. The
spinor bundles are defined to be 
\begin{equation*}
\begin{split}
  S^+_{\mathfrak{c}} := & \Lambda^{0,0}(X) \oplus \Lambda^{0,2}(X) \ ,\\
  S^-_{\mathfrak{c}} := & \Lambda^{0,1}(X) \ ,
\end{split}
\end{equation*}
and the Clifford multiplication is given by 
\begin{equation*}
\begin{split}
 \gamma  & :  \Lambda^1(T^*X) \to
\Hom_{\C}(S^+_{\mathfrak{c}},S^-_{\mathfrak{c}}) \\
   & u \mapsto \sqrt{2} (e(u^{0,1}) - e^*(u^{0,1}))  \ .
\end{split}
\end{equation*}
The induced isomorphism
\[
\gamma : \Lambda^2_+(X) \tensor \C \to \mathfrak{sl}(S^+_{\mathfrak{c}})
\]
is then seen to be given by the formula
\begin{equation}\label{cliffordmultiplication}
\gamma(\eta^{1,1} + \eta^{2,0} + \eta^{0,2}) = 
4 \begin{pmatrix} -i \Lambda_g(\eta^{1,1}) & - * (\eta^{2,0} \wedge \_ \ )
\\ \eta^{0,2} & i \Lambda_g(\eta^{1,1}) \end{pmatrix} \ .
\end{equation}
Here we use the commonly used convention to denote contraction with $\omega_g$,
that is $e^*(\omega_g)$, by the symbol $\Lambda_g$.

Now suppose that $X$ is a K\"ahler surface. This means that first the almost
complex structure $J$ is integrable to a complex structure, and second that the
K\"ahler form $\omega_g$ is closed, $d \omega_g = 0$. The condition of
closedness implies (cf. \cite{KN2}, p. 148) that the the almost complex
structure $J$ is parallel with respect to the Levi-Civita-connection
$\nabla_g$. As a consequence, the splittings 
\[
\Lambda^k(X)\tensor \C =
\oplus_{p+q=k}\Lambda^{p,q}(X)
\]
are $\nabla_g$-parallel, where we also denote by $\nabla_g$ the connection induced by the Levi-Civita
connection on all exteriour powers of $T^*X$. The canonical $Spin^c$-connection is now simply given by the
the connection $\nabla_g$ in the bundles $\Lambda^{0,0},
\Lambda^{0,1} $ and $\Lambda^{0,2}$.

%

Let $E$ be a Hermitian vector bundle on $X$, and further $\nabla_{\hat{A}}$
a unitary connection on $E$. We shall use the notation convention
$\Lambda^{p,q}(E) := \Lambda^{p,q}(X) \tensor E$, and by $\Omega^{p,q}(E)$ we
shall denote the space of sections of the latter bundle,  $\Omega^{p,q}(E) =
\Gamma(\Lambda^{p,q}(E))$ .
\begin{definition}
  The operator $\bar{\partial}_{\hat{A}} : \Omega^{p,q}( E)
\to \Omega^{p,q+1}( E)$ is defined to be the composition of 
$d_{\hat{A}} : \Omega^{p+q}(E) \to
\Omega^{p+q+1}(E)$, the extension of the exteriour derivative to forms with values in $E$ by means of the connection $\nabla_{\hata}$,
 with the bundle projection
$\Lambda^{p+q+1}(E) \to
\Lambda^{p,q+1}( E)$.
\end{definition}

The Dirac operator associated to the canonical
$Spin^c$-connection
$\nabla_g$ in the canonical $Spin^c$-structure $\mathfrak{s}_c$ and the unitary
connection $\hat{A}$ in the Hermitian bundle $E$ 
is expressible in terms of the above operator $\bar{\partial}_{\hat{A}}$ and
its formal $L^2$-adjoint $\bar{\partial}_{\hat{A}}^*$ as follows:
\begin{equation} \label{dirackaehler}
  D_{\hat{A}} = \sqrt{2} \left( \bar{\partial}_{\hat{A}} +
\bar{\partial}_{\hat{A}}^* \right) \ .
\end{equation}
This is a well-known fact in the case $n=1$ \cite{Hi}. The proof of the general case follows along the same
lines. In particular, the proof given in the lecture notes \cite{T3} is directly applicable to our
situation. \qed 

%

We will now study the $U(n)$ monopoles associated to the data
$(\mathfrak{c},E)$ with spinor bundles $W_{\mathfrak{c},E}^\pm = S^\pm \tensor
E$. Note that, up to tensoring $E$ with a line bundle, we can always
assume that general data $(\mathfrak{s},E)$ is of the particular form
$(\mathfrak{c},E)$. Now according to the isomorphism
$W^+_{\mathfrak{c},E} \cong \Lambda^{0,0}(E) \oplus \Lambda^{0,2}(E)$ a
spinor $\Psi \in \Gamma(X;W^+_{\mathfrak{c},E})$ can be written as
$\Psi=(\alpha,\beta)$ with $\alpha \in \Omega^{0,0}(X;E)$ a section of $E$ and
$\beta \in \Omega^{0,2}(X;E)$ a 2-form of type $(0,2)$ with values in $E$.
We introduce the following notations. We denote by $ ^{-}:
\Lambda^{p,q}(E) \to \Lambda^{q,p}(E^*)$ the conjugate linear isomorphism which
is the tensor product of complex conjugation on the forms and the conjugate
linear isomorphism specified by the hermitian structure on the bundle $E$. We
denote by $^* : \Lambda^{p,q}(E) \to \Hom(\Lambda^{p,q}(E),\C)$ the conjugate
linear isomorphism specified by the Hermitian structure on $\Lambda^{p,q}(E)$.
For an endomorphism $f \in \End(E)$ we denote $\{ f \}_\tau := (f)_0 +
\frac{\tau}{n} \tr(f) \id_E $, where $(f)_0$ denotes the trace-free part of
$f$. Thus we simply have $\{f\}_1 = f$. With this said we can write
$\mu_{0,\tau}(\Psi)$ according to the above isomorphism as 

\begin{equation}\label{mukaehler}
  \mu_{0,\tau}(\Psi) = \begin{pmatrix}
  \frac{1}{2} \left( \{\alpha \alpha^*\}_\tau - \{ * \beta \wedge \bar{\beta}
\}_\tau \right)  & \{\alpha \beta^*\}_\tau \\
\{\beta \alpha^*\}_\tau & \frac{1}{2} \left( \{\beta \beta^*\}_\tau - \{ \alpha
\alpha^*\}_\tau \right)
              \end{pmatrix} \ .
\end{equation}
It is worth pointing out here that we have $\beta \beta^* = * \beta \wedge
\bar{\beta}$ which is true because $\Lambda^{0,2}$ is 1-dimensional. In other
words, the two diagonal entries only ``look" differently. 

With the above formulae (\ref{cliffordmultiplication}) we can now
write down the monopole equations (\ref{monopole equations}) with parameter
$\tau$ and as perturbation the imaginary-valued self-dual 2-form $\eta$ for the pair consisting of the spinor
$\Psi = (\alpha,\beta) \in \Gamma(X;\Lambda^{0,0}(E) \oplus \Lambda^{0,2}(E))$
and the connection $\hat{A}$ in $E$:

\begin{equation}\label{un equations kaehler} 
\begin{split}
  	\bar{\partial}_{\hat{A}} \alpha + \bar{\partial}^*_{\hat{A}} \beta & = 0
\\
	F_{\hat{A}}^{0,2} & =  \frac{1}{4} \{\beta \alpha^*\}_\tau + 4
	\eta^{0,2}
\\
  	- i \Lambda_g (F_{\hat{A}}) & =  
	\frac{1}{8} \left\{ \alpha \alpha^* - *(\beta \wedge
	\bar{\beta})\right\}_\tau \	- i \Lambda_g(\eta)  
	\ .
\end{split}
\end{equation}
Indeed, the curvature equation of (\ref{monopole equations}) splits into four
equations according to the above splitting, but the two equations resulting
from the diagonal entries are equivalent, and, using that
$\overline{F_{\hat{A}}^{0,2}} = -F_{\hat{A}}^{2,0}$ (here again, $^-$ denotes
the complex-conjugation on the forms and the hermitian adjoint on $\End(E)$),
the two off-diagonal equations also prove to be equivalent.

\subsection{Decoupling phenomena, moduli spaces for $b_2^+(X) >  1$ and holomorphic 2-forms}
As mentioned before a lot of the analysis of the classical monopole equations on K\"ahler surfaces carries over to our situation. Before we consider the perturbed monopole equations we shall first draw some intermediate conclusions from the unperturbed monopole equations.
In particular there is a decoupling result completely analogous to the classical situation, interpreting monopoles as `vortices', c.f. also \cite{BG}, \cite{T}. 

\begin{prop}\label{decoupling}
Let $X$ be a K\"ahler surface. Suppose that the configuration $(\Psi,\hat{A})
\in \Gamma(X;S^+_\mathfrak{c} \tensor E) \times \mathscr{A}(E)$ solves the {\em
unperturbed} $U(n)$ monopole equations with parameter $\tau \in [0,1]$. If we
write the spinor as $\Psi = (\alpha,\beta)$ according to the decomposition
$S^+_\mathfrak{c} \tensor E \cong \Lambda^{0,0}(E) \oplus \Lambda^{0,2}(E)$ then
one of the following two statements holds:
\begin{enumerate}
\item The second factor of the spinor vanishes identically, $\beta \equiv 0$.
Furthermore the pair  $(\alpha,\hat{A})$ satisfies the following
`Vortex-type' equations
\begin{equation}\label{vortexalpha}
\begin{split}
  \bar{\partial}_{\hat{A}} \alpha & = 0 \\
  F_{\hat{A}}^{0,2} & = 0 \\
  i \Lambda_g (F_\hata) & = - \frac{1}{8} \{ \alpha \alpha^*\}_\tau \ .
\end{split}
\end{equation}
\item The first factor of the spinor vanishes identically, $\alpha \equiv 0$.
Furthermore the pair $(\beta,\hat{A})$ satisfies the following equations
\begin{equation}\label{vortexbeta}
\begin{split}
  \bar{\partial}_{\hat{A}}^* \beta & = 0 \\
  F_{\hat{A}}^{0,2} & = 0 \\
  i \Lambda_g (F_\hata) & = + \frac{1}{8} \{\beta \beta^*\}_\tau \ .
\end{split}
\end{equation}
\end{enumerate}
\end{prop}
{\em Proof:} 
Using the first two of the monopole equations (\ref{un
equations kaehler}) we get:
\begin{equation*}
\dbar_\hata \dbar_\hata^* \beta  = - \dbar_\hata \dbar_\hata \alpha =
-F_\hata^{0,2} \alpha = -\frac{1}{4} \{\beta \alpha^*\}_\tau \alpha \ .
\end{equation*}
We take the inner product with $\beta$ to get now:
\begin{equation*}
\begin{split}
 \klammer{\beta,\dbar_\hata \dbar_\hata^* \beta} & = -\frac{1}{4}
\klammer{\beta, \geklatau{\beta \alpha^*} \alpha} \\
& = -\frac{1}{4} \klammer{\abs{\beta}^2 \abs{\alpha}^2 - \frac{1-\tau}{n}
\klammer{\beta, \tr(\beta \alpha^*) \alpha}} \\
& \leq \klammer{-\frac{1}{4} + \frac{1-\tau}{4n}} \abs{\alpha}^2 \abs{\beta}^2
\\
& \leq 0 \ .
\end{split}
\end{equation*}
Here we have used the Cauchy-Schwarz inequality, noting also that
$\abs{\tr(\beta \alpha^*)} \leq \abs{\beta \alpha^*} = \abs{\beta}
\abs{\alpha}$. Integrating now the latter inequality over the whole manifold $X$
yields the following:
\begin{equation}\label{conclusion inequality}
 0 \ \leq \norm{\dbar_\hata^* \beta}^2 \ \leq \klammer{-\frac{1}{4} +
	\frac{1-\tau}{4n}} \int_X \abs{\alpha}^2 \abs{\beta}^2 vol_g \leq 0
\end{equation}
Thus we get $\dbar_\hata^* \beta = 0$ and from the Dirac equation also
$\dbar_\hata \alpha = 0$. If further we have $\tau > 1-n$ then we see from the
last inequality that at any point of the manifold $X$ we have  $\alpha = 0$ or
$\beta = 0$. But we have $ 0 = \dbar_\hata^* \dbar_\hata \alpha =
\Delta_{\dbar_\hata} \alpha$, and because $\Delta_{\dbar_\hata}$ is an elliptic
second order operator with scalar symbol it follows from Aronaszajin's theorem
\cite{A} that solutions to $\Delta_{\dbar_\hata} \alpha = 0$ satisfy a unique
continuation theorem. Similarly we have $ 0 = \dbar_\hata \dbar_\hata^* \beta =
\Delta_{\dbar_\hata} \beta$, so the same holds for $\beta$. Therefore, if one
of $\alpha$ or $\beta$ vanishes on an open subset of $X$, then it vanishes on
the whole of $X$. The conclusions now follow from (\ref{un equations kaehler}).
\qed
\begin{remark}
If $\tau \neq 0$ the moduli space $M_{\mathfrak{c},E}(\tau,0)$
can only contain either solution with $\alpha \neq 0$ or with $\beta \neq 0$. 
This follows from taking the trace of the third equation of
(\ref{un equations kaehler}) and then integrating it over the whole manifold.
The left hand term yields then the topological quantity $- 2 \pi \, \langle
c_1(E) \smile \eckklammer{\omega_g},[X] \rangle$. 
\end{remark}

On a K\"ahler surface we have $\Delta = 2
\Delta_{\dbar}$, just reflecting again the compatibility between the complex structure and the Riemannian
metric. Therefore the harmonic differential forms are also $\dbar$-harmonic and vice versa. In
particular, we get the following decomposition from the Hodge-theorem:
\begin{equation}\label{h2_kaehler}
H^2_{dR}(X;\C) = H^{2,0}_{\bar{\partial}}(X) \oplus H^{1,1}_{\bar{\partial}}(X)
\oplus H^{0,2}_{\bar{\partial}}(X) \ .
\end{equation}

\begin{corollary}
If there are solutions to the unperturbed $U(n)$-monopole equations associated
to the data $(\mathfrak{c},E)$ and to the parameter $\tau \in [0,1]$, then the
image $c_1^\R(E)$ in real (complex) cohomology of the first Chern-class $c_1(E)
\in H^2(X;\Z)$ is of type $(1,1)$ according to the above decomposition
(\ref{h2_kaehler}).
\end{corollary}
{\em Proof:}
Under these conditions there is a connection $\hat{A}$ on $E$ with
$F_\hata^{0,2} = 0 = F_\hata^{2,0}$. From the Chern-Weil formula we have that 
$\frac{-1}{2 \pi i} \eckklammer{\tr (F_\hata)} = c_1^\R(E)$. There is a 1-form
$\lambda$ such that $\omega:=\tr(F_\hata) - \dbar \lambda$ is $\dbar$-harmonic
and this class also represents $c_1^\R(E)$. We have $\omega^{2,0} = 0$ and
$\omega^{0,2} = \dbar \lambda$. But a class is $\dbar$ - harmonic if and only 
each component according to $\Omega^{p+q}(X) = \oplus \Omega^{p,q}(X)$ is
$\dbar$ - harmonic. But then the harmonic form $\omega^{2,0} = \dbar \lambda$
must be zero, as it is a $\dbar$ - exact form also. \qed

In the classical theory a common perturbation of the monopole equations was to
perturb with imaginary-valued self-dual 2-forms $\eta$ such that $\eta^{2,0}$ is
a holomorphic form \cite{W} \cite{Bq}. There are such forms with $\eta^{2,0}
\neq 0$ precisely if $b_2^+(X) > 1$. We will now consider this type of
perturbation in the general case of $U(n)$ monopoles even though these perturbations are not enough to get generic regularity of the moduli space in the case $n > 1$.  However, it will turn out that the moduli spaces perturbed in this way are empty in the case $n > 1$ as soon as the perturbing form $\eta$  is non-zero. 

If the unperturbed $U(n)$ monopole moduli space is empty then any invariant 
derived by the scheme `evaluation of cohomology classes on the fundamental cycle of the moduli space' should be zero. Indeed, that kind of invariant would be defined with a `generic' moduli space, i.e. one which is cut out transversally by the suitably perturbed monopole equations. An empty moduli space is always generic. Thus if there is a non-trivial invariant derived from some generic moduli space then the associated unperturbed moduli space may be not generic, but it could not be empty. Therefore it is natural to consider topological data $(\data)$ only for situations where the unperturbed $U(n)$ monopole moduli spaces are a priori non-empty. As we have seen, this can only be the case if the first Chern-class $c_1^\R(E)$ is of type $(1,1)$ according to the decomposition (\ref{h2_kaehler}). Therefore we shall include this hypothesis to the next two results, the following theorem and its corollary:

\begin{theorem}
Let $X$ be a K\"ahler surface and let $E$ be a bundle such that its first
Chern-class $c_1^\R(E)$ is of type $(1,1)$. Further let $\eta$ be an
imaginary-valued 2-form with $\eta^{2,0}$ holomorphic. Then the $U(n)-$
monopole equations (\ref{un equations kaehler}) associated to the data
$(\mathfrak{c},E)$, to the perturbation form $\eta$, and to the parameter $\tau \in (0,1]$ are equivalent to the
following system of equations:
\begin{equation}\label{equations kaehler split}
\begin{split}
  	\dbar_{\hat{A}} \alpha & = 0 \\
	\dbar_{\hat{A}}^* \beta & = 0 \\
	F_{\hat{A}}^{0,2} & = 0 \\
	\frac{1}{4} \{\beta \alpha^*\}_\tau & = \eta^{0,2} \\
	-i \Lambda_g(F_{\hat{A}}) & = \frac{1}{8} \{ \alpha \alpha^* - \beta
	\beta^* \}_\tau - i \Lambda_g(\eta) \\
\end{split}
\end{equation}
\end{theorem}
{\em Proof: } 
We will derive the following formula for a
solution $((\alpha,\beta),\hata))$ to the $U(n)-$monopole equations (\ref{un
equations kaehler}) with parameter $\tau$ and perturbation $\eta$:
\begin{equation}\label{formula}
\begin{split}
0 =  & 4 \, \norm{F_\hata^{0,2}}^2_{L^2(X)} \ 
+ 4 \, \frac{1-\tau}{\tau n} \, \norm{\tr F_\hata^{0,2}}^2_{L^2(X)}
+ \norm{\dbar_\hata^*
\beta}^2_{L^2(X)} \\
	& - \frac{4}{\tau} \, \langle 2 \pi i \ [\eta^{2,0}] \smile c_1(E) , [X]
	\rangle \
	 \ .
\end{split}
\end{equation}
The conclusion then clearly follows as the topological term vanishes by assumption.
\\

Provided that we have $\tau \neq 0$ the endomorphism $\beta \alpha^*$ can be
expressed as
\begin{equation}\label{expression of ab}
\begin{split}
\beta \alpha^* & = \geklatau{\beta \alpha^*} \ + \frac{1-\tau}{n} \tr (\beta
		\alpha^*) \\
		& = \geklatau{\beta \alpha^*} \ + \frac{1-\tau}{\tau n} \tr
		(\geklatau{\beta \alpha^*}) \\
		& = 4 \, F_\hata^{0,2} - 4\,  \eta^{0,2} + 4 \, \frac{1-\tau}{n \tau} \,
 			\tr(F_\hata^{0,2}) - 4 \, \frac{1-\tau}{\tau} \, \eta^{0,2} 	\ ,
\end{split}
\end{equation}
where the last equation used the second of the monopole equations (\ref{un
equations kaehler}) and the trace of it. 

Again we get from the Dirac-equation that $
\dbar_\hata \dbar_\hata^* \beta + F_\hata^{0,2} \alpha = 0 $, 
so that after taking the pointwise inner-product with $\beta$ and using the
above equation (\ref{expression of ab}) we get: 
\begin{equation}\label{to integrate}
\begin{split}
0 &  = \klammer{\beta, F_\hata^{0,2} \alpha} \ + \klammer{\beta, \dbar_\hata
		\dbar_\hata^* \beta } \\
	& = \klammer{ \beta \alpha^* , F_\hata^{0,2} } \ + \klammer{\beta, \dbar_\hata \dbar_\hata^* \beta } \\ 
	& = 4 \, \abs{F_\hata^{0,2}}^2 \ - 4 \, \klammer{\eta^{0,2}, F_\hata^{0,2} } \\
	& \ \ \ \ \ 	+ 4 \frac{1-\tau}{n \tau} \abs{ \tr(F_\hata^{0,2})}^2 \ - 4  \, \frac{1-\tau}{\tau} \, 
		\klammer{\eta^{0,2}, F_\hata^{0,2}} + \,  \klammer{\beta, \dbar_\hata
		\dbar_\hata^* \beta } \\
	& = 4  \, \abs{F_\hata^{0,2}}^2 \  + 4 \, \frac{1-\tau}{n \tau} \abs{ \tr(F_\hata^{0,2})}^2  -
\frac{4}{\tau} \klammer{\eta^{0,2}, F_\hata^{0,2} } + \,  \klammer{\beta, \dbar_\hata
		\dbar_\hata^* \beta }
\end{split}
\end{equation}
As the next step we will integrate this whole equation over $X$. Beforehand we shall remark that  $\eta^{2,0}$ is closed, and therefore the following integral is of topological nature:
\begin{equation}\label{topological term}
\begin{split}
\int_X \klammer{\eta^{0,2}, F_\hata^{0,2} } \, vol_g & = \int_X \bar{\eta^{0,2}} \wedge * \tr(F_\hata^{0,2}) \\
	& = - \, \int_X \eta^{2,0} \wedge \tr(F_\hata^{0,2}) \\
	& = \, 2 \pi i \ \langle  [\eta^{2,0}] \smile c_1(E) , [X]
	\rangle
\end{split}
\end{equation}
With this said the integral of the formula (\ref{to integrate}) clearly yields the above formula (\ref{formula}).

\qed

\begin{corollary}\label{final conclusion}
Let $X$ be a K\"ahler surface with $b_2^+(X) > 1$ and let $E$ be a bundle such
that its first Chern-class $c_1^\R(E)$ is of type $(1,1)$. Then for any
self-dual imaginary valued 2-form $\eta$ with $\eta^{2,0}$ holomorphic and {\em
non-zero} and constant $\tau \in (0,1]$ the moduli space $M_{\mathfrak{c},E}(\eta,\tau)$ is empty. 
\end{corollary}
{\em Proof:}
  Under the given hypothesis the preceeding theorem implies that 
\begin{equation}
\label{impossible}
\geklatau{ \beta \alpha^*} = 4 \eta^{0,2} \ \id_E \ .
\end{equation}
But using the definition of $\geklatau{\beta \alpha^*}$ it is a pure matter
of linear algebra to check that for $\eta^{0,2} \neq 0$ this is impossible if $n
\geq 2$, because the left hand side of the equation (\ref{impossible}) can
never be a mutliple of the identity, unless $\alpha = 0$ or $\beta = 0$.

\section{Appendix}
{\em Proof of Proposition \ref{properness}:} Obviously we have
$\abs{\mu_{0,\tau}(\Psi)} \geq \abs{\mu_{0,0}(\Psi)}$, so for the first
assertion it will be enough to consider $\mu_{0,0}$ alone. We will show that
$\mu_{0,0}(\Psi) =
0$ implies $\Psi = 0$. Because $\mu_{0,0}$ is quadratic and the unit sphere
inside
$\C^2 \tensor \C^n$ is compact we then get the claimed uniform
properness-inequality (\ref{propernessconstant}). 

We shall use the canonical
isomorphism $\C^2
\tensor \C^n \cong \C^n \oplus \C^n$, which permits to write a general element 
\[
  \Psi=\begin{pmatrix} 1 \cr 0 \end{pmatrix} \tensor \alpha + 
      \begin{pmatrix} 0 \cr 1 \end{pmatrix} \tensor \beta 
\]
as 
\[
  \Psi= \begin{pmatrix}\alpha \cr \beta \end{pmatrix} \ .
\]
We then have
\begin{equation*}
\begin{split}
\mu_{0,0}(\Psi) & =  P (\Psi \Psi^*) \\
	& =  P \left( \begin{pmatrix} \alpha \cr \beta \end{pmatrix}
\begin{pmatrix} \alpha^* & \beta^*\end{pmatrix} \right) \\
	& =  P \begin{pmatrix} \alpha \alpha^* & \alpha \beta^* \cr 
			\beta \alpha^* & \beta \beta^* \end{pmatrix} \\
	& =   \begin{pmatrix} \frac{1}{2} (\alpha \alpha^* - \beta \beta^*)_0 & 
	      (\alpha \beta^*)_0 \cr (\beta \alpha^*)_0 & \frac{1}{2}(\beta
\beta^* - \alpha \alpha^*)_0 \end{pmatrix} \ . 
\end{split}
\end{equation*}
In particular, if $\mu_{0,0}(\Psi)= 0$, then we have $(\alpha \beta^*)_0 = 0$. 
\begin{lemma} \label{lemmainlemma}
The equation $(\alpha \beta^*)_0 = 0$ implies that $\alpha=0$ or $\beta=0$. In
other words, the bilinear map $(\alpha, \beta) \to (\alpha \beta^*)_0$ is
without
zero-divisors (here $n \geq 2$ is implicitely understood).
\end{lemma}
{\em Proof of Lemma \ref{lemmainlemma}:} Write the elements $\alpha$ and $\beta
$ as 
\[
  \alpha = (\alpha_i)_{i=1}^n \ , \quad \beta= (\beta_i)_{i=1}^n \ .
\]
Then the equation $(\alpha \beta^*)_0 = 0 $ reads in  matrix-form
\begin{equation*}
\begin{split}
\left( \begin{array}{ccc} \alpha_1 \bar{\beta_1} - \frac{1}{n} \sum \alpha_i
\bar{\beta_i} & \dots & \alpha_1 \bar{\beta_n} \\
\vdots & \ddots & \vdots \\
\alpha_n \bar{\beta_1} & \dots & \alpha_n \bar{\beta_n}-\frac{1}{n} \sum
\alpha_i \bar{\beta_i} \end{array} \right) = 0  \ .
\end{split}
\end{equation*}
Suppose $\beta\neq 0$, for instance $\beta_j \neq 0$. Then the $j^{\text{th}}$
column implies that $\alpha_i= 0$ for all $i \neq j$. Thus the
$j^{\text{th}}$ element in the $j^{\text{th}}$ column simplifies,
\begin{equation*}
 \alpha_j \bar{\beta_j} - \frac{1}{n}\sum_{i=1}^{n} \alpha_i \bar{\beta_i} \, =
\, (1-\frac{1}{n}) \alpha_j \bar{\beta_j} \ .
\end{equation*}
Therefore we have $\alpha_j= 0$ as well, so that we have $\alpha=0$. The case
$\alpha \neq 0 $
is analogous. \qed
Returning to the problem
\[
  \begin{pmatrix} \frac{1}{2} (\alpha \alpha^* - \beta \beta^*)_0 & 
	      (\alpha \beta^*)_0 \cr (\beta \alpha^*)_0 & \frac{1}{2}(\beta
\beta^* - \alpha \alpha^*)_0 \end{pmatrix} = 0 \ ,
\]
we see that the lemma gives $\alpha = 0$ or $\beta = 0$. Suppose, without loss
of
generality, that the first is the case. Then we are left with $(\beta \beta^*)_0
= 0 $. Now again with lemma \ref{lemmainlemma} we see that this also implies
$\beta=
0$. Therefore 
\[
\Psi= \begin{pmatrix} \alpha \cr \beta \end{pmatrix} = 0 \ .
\]

The second assertion now follows from the first, remembering that $P$ and $Q$
are both orthogonal projections. For non-negative $\tau$ we have the
inequality 
\begin{equation*}
\begin{split}
\left(\mu_{0,\tau}(\Psi)\Psi,\Psi \right) = & \, \left(P (\Psi \Psi^*) \Psi,
\Psi
\right) + \tau \klammer{Q(\Psi \Psi^*) \Psi, \Psi} \\
  = & \, \left( P (\Psi \Psi^*), \Psi \Psi^* \right) + \tau \klammer{Q(\Psi 
\Psi^*), \Psi \Psi^*} \\
  = & \, \left( P (\Psi \Psi^*), P (\Psi \Psi^*) \right) + \tau \klammer{Q(\Psi
\Psi^*), Q(\Psi \Psi^*)} \\
  \geq & \ \abs{\mu_{0,0}(\Psi)}^2 \\
  \geq & \ c^2 \abs{\Psi}^4 \ .
\end{split}
\end{equation*}
\qed

\end{document}